\documentclass[a4paper,11pt,reqno]{amsart}

\begin{document}

\newcommand{\col}[1]{{\left(\begin{smallmatrix}{#1}\\1\end{smallmatrix}\right)}}
\newcommand{\colinf}{{\left(\begin{smallmatrix}{1}\\0\end{smallmatrix}\right)}}
\newcommand{\mat}[4]{%
{\left(\begin{smallmatrix}{#1}&{#2}\\{#3}&{#4}\end{smallmatrix}\right)}}
\newcommand{\tr}{{\rm tr}}
\newcommand{\rat}[3]{{\left(\frac{{#1}{#2}+{#1}+{#3}}{{#1}+{#2}}\right)}}

\author{Andrew Thomason}
\address{Department of Pure Mathematics and Mathematical Statistics,
Centre for Mathematical Sciences,
Wilberforce Road, Cambridge CB3 0WB, England}

\email{a.g.thomason@dpmms.cam.ac.uk}

\date{6th March 2015}


\title{A Paley-like graph in characteristic two}

\begin{abstract}
  The Paley graph is a well-known self-complementary pseudo-random graph,
  defined over a finite field of odd order. We describe an attempt at an
  analogous construction using fields of even order. Some properties of the
  graph are noted, such as the existence of a Hamiltonian decomposition.
\end{abstract}

\maketitle

\section{Introduction}

The well-known Paley graph is a pseudo-random graph whose vertex set is the
finite field $F_q=GF(q)$ of order $q\equiv1$~(mod 4). The pair of vertices
$a$, $b$ is joined by an edge if $a-b$ is a square in $F_q$. Since $-1$ is
a square the graph is well defined. It follows from elementary properties
of the quadratic character that the Paley graph is vertex-transitive,
self-complementary, and each edge is in $(q-5)/4$ triangles. A graph of
order $q$ in which every edge is in $(q-5)/4$ triangles and whose
complement has the same property is sometimes called a {\em conference}
graph. The Paley graph is thus {\em ipso facto} a pseudo-random graph, as
explained in detail in~\cite{T}, and in a somewhat less quantitative
fashion in Chung, Graham and Wilson in~\cite{CGW}.

The other odd prime powers, namely those where $q\equiv3$~(mod 4), cannot
be used to construct Paley graphs since $-1$ is not a square. However this
very property allows the construction of a {\em tournament}, or oriented
complete graph, on the vertex $F_q$ by inserting an edge oriented from $a$
to $b$ if $a-b$ is a square. Since $-1$ is not a square, exactly one of
$a-b$ and $b-a$ is a square, so we do indeed construct a tournament (Graham
and Spencer~\cite{GS}).

The property of pseudo-randomness, even when quantified, does not suffice to
give all the information that one would like to have about the Paley
graphs; in particular, it is not known what the clique number is. When $q$
is prime the calculation reduces to difficult and so far unsolved problems
involving the estimation of character sums (though if $q$ is a square the
clique size is exactly $\sqrt q$). For more information,
see~\cite[Section~2.5.1]{TR}.

As regards Hamiltonian cycles, the Paley graphs hold fewer secrets. If $q$ is
a prime then the Paley graph is a circulant graph and, since the edges of a
given distance in a circulant of prime order form a Hamiltonian cycle, it
follows that the Paley graph in this case is not just Hamiltonian but it
has a {\em Hamiltonian decomposition}, that is, its edge set is the union
of edge-disjoint Hamiltonian cycles.

In a finite field of characteristic two, every element is a square, and the
definition of the Paley graph is of little value. From the graph
theoretical point of view, though, characteristic two has some innate
attraction. In this note we describe an attempt to find different graphs,
similar in spirit to the Paley graphs but defined in relation to the field
$F_q$ for even~$q$, which are vertex-transitive and self-complementary. We
might even hope to find a graph which is a conference graph, or which is
more easily analysed than the Paley graph. Whilst these more ambitious aims
are not realised, we do describe the more accessible properties of the graphs.

\section{Definition}\label{defn}

We begin with a definition of the graphs, and defer to~\S\ref{hist} a
discussion of what lies behind it. Choose as vertex set $V$ the
elements of $PG(1,q)$, the projective space of dimension one over $F_q$. We
label the elements of $V$ in the natural way, namely
$$V=\left\{\colinf,\col0,\col1,\ldots,\col x,\ldots\right\}
=\{\infty,0,1,\ldots,x,\ldots\}.$$

Given an element $x\in F_q$ its {\em trace} is defined to be
$\tr(x)=x+x^2+x^4+\ldots+x^{q/2}$.  Let $q=2^k$ and let $a$ be an element of
$F_q$ with $\tr(a)=1$. For even $k$ we define a graph $G_k(a)$ on the vertex
set $V$ by 
$$xy\in E(G_k(a))\qquad \mbox{if}\qquad\tr\rat{x}{y}{a}=0.$$
For odd $k$ we define a tournament $G_k(a)$, having an edge directed from $x$
to $y$ whenever the same equation is satisfied.

We shall show in~\S\ref{tracing} that $G_k(a)$ is well defined. Moreover,
although $G_k(a)$ as a labelled graph depends on the value of~$a$ (for example,
the neighbourhood of the vertex~0 is the set of elements $y$ such that
$\tr(a/y)=0$), we shall show in~\S\ref{elem} that all the graphs (or
tournaments) so defined are isomorphic. This allows us to refer to any member
of this collection of graphs as {\em the graph $G_k$} when there is no danger
of confusion.

At first appearance the definition of the graph $G_k$ looks somewhat
contrived. We attempt in~\S\ref{hist} to show that the definition does in
fact arise fairly naturally. Having made a few elementary remarks about the
properties of~$F_q$ (in~\S\ref{field}) we establish in~\S\ref{elem} that
$G_k(a)$ is a vertex-transitive self-complementary graph whose isomorphism
class is independent of~$a$, as claimed. Finally, we explore some of the
further properties of the graph $G_k$; in particular we show that it is a
pseudo-random graph, having a Hamiltonian decomposition.

\section{Background}\label{hist}

The Paley graph is a circulant graph when $q$ is prime; that is, its
vertices may be labelled $\{0,1,\dots,q-1\}$, and whether $xy$ is an edge
depends only on the difference $|x-y|$. It is therefore necessarily
vertex-transitive, and it is also self-complementary. A regular
self-complementary graph has order $\equiv 1$ (mod~4), and obviously there
is no such graph with vertex set $F_q$ when $q$ is even. We set out to
define a circulant graph on vertex set $PG(1,q)$.

The group $PSL(2,q)$, comprising the $2\times2$ matrices of determinant
one, acts on $V=PG(1,q)$. As usual, we associate with the matrix $\mat
abcd$ the M\"obius, or linear fractional, map $z\mapsto (az+b)/(cz+d)$.  We
use two simple facts about these maps; that they form a group, and (for
this background discussion) that a map is determined by its action on any
three points. For convenience and completeness we assume only a minimal
familiarity with properties of finite fields. Much more can be found in the
classical algebraic text of Dickson~\cite{D} or the more recent and
geometrical Hirschfeld~\cite{H}. Both these authors pay attention to
the characteristic two case needed here.

In order to begin constructing a circulant on $V$ we need a M\"obius
transformation of order $q+1$. It is not hard to show, though we don't need
this fact, that every transformation with no fixed points is conjugate to
one of the form $z\mapsto a/(z+1)$ such that the equation $x^2+x=a$ has no
solution in $F_q$. Let us then consider such a map. The condition that
$x^2+x=a$ has no solution is equivalent to the condition $\tr(a)=1$
(see~\S\ref{tracing}). Amongst such transformations there exist some of order
$q+1$ (see~\S\ref{conjugal}).

Take such a transformation $\alpha$. Then $V=\{\infty,\alpha(\infty),
\alpha^2(\infty),\ldots, \alpha^q(\infty)\}$. For convenience, we write
$v_i=\alpha^i(\infty)$, so $V=\{v_0,v_1,\ldots,v_q\}$.  Notice that, for
example, $v_1=\alpha(\infty)=0$ and $v_2=\alpha(0)=a$. Moreover
$\alpha^{-1}(z)=1+a/z$, so $v_q=\alpha^q(\infty)=\alpha^{-1}(\infty)=1$ and
$v_{q-1}=\alpha^{-1}(1)=1+a$. It is easily verified, by induction on~$i$,
that $v_{q-i}=1+v_i$ (the induction step being
$v_{q-i-1}=\alpha^{-1}(v_{q-i}) =1+a/v_{q-i} = 1+a/(1+v_i) = 1+\alpha(v_i)=
1+v_{i+1}$). Subscripts may be reduced modulo $(q+1)$, so we write
$v_{-i}=1+v_{i}$.

We may, therefore, define a circulant graph on~$V$ as follows. Choose a map
$f:F_q\to F_2$, to be specified later. The neighbours of $\infty=v_0$ will
be those $v_i$ for which $f(v_i)=0$. In general, $v_iv_j$ will be an edge
if $v_0v_{j-i}$ is an edge, which is to say, if $f(v_{j-i})=0$. In order
that the graph be well defined we must ensure that $f(v_{i-j})=f(v_{j-i})$,
which we have seen is equivalent to $f(x+1)=f(x)$. (This section is just to
motivate the earlier definition, so we ignore tournaments here.)

Let us see how to compute whether $xy$ is an edge, given $x,y\in V$. Let
$x=v_i$ and $y=v_j$. Then $xy$ will be an edge if
$f(v_{i-j})=f(v_{j-i})=0$. Now $v_{i-j}=\alpha^{-j}(x)$. We claim that the
map $\alpha^{-j}$ is identical to the M\"obius map
$\beta(z)=(zy+z+a)/(z+y)$, and so $v_{i-j}=\beta(x)=(xy+x+a)/(x+y)$.  To
check the claim, it suffices to show that $\alpha^{-j}$ and $\beta$ act
identically on the three distinct points $v_{j-1}$, $v_j$ and
$v_{j+1}$. Now $\alpha^{-j}$ maps these points to $v_{-1}=1$, $v_0=\infty$
and $v_1=0$. But $v_j=y$, $v_{j-1}=\alpha^{-1}(y)=1+a/y$ and
$v_{j+1}=\alpha(y)=a/(y+1)$. Thus $\beta(v_{j-1})=\beta(1+a/y)=1$,
$\beta(v_j)=\beta(y)=\infty$ and $\beta(v_{j+1})=\beta(a/(y+1))=0$, proving
the claim. We conclude that $xy$ is an edge if
$f((xy+x+a)/(x+y))=f(v_{j-i})=0$.

The map $v_i\mapsto v_{2i}$ is a permutation of~$V$ which leaves $v_0$
fixed. An easy way to ensure that our circulant graph is self-complementary
is to arrange that this map interchanges the graph with its complement. So
we wish to arrange that if $x=v_i$ then $f(v_{2i})\ne f(x)$, or,
equivalently, $f(v_{2i})+f(x)=1$. If we put $j=-i$ then $v_{2i}=v_{i-j}$,
and using the calculation in the previous paragraph with $y=v_{-i}=1+x$, we
see that $v_{2j}=(xy+x+a)/(x+y)=x^2+a$.

Therefore this procedure will yield a self-complementary vertex-transitive
graph if we select a function $f:F_q\to F_2$ such that $f(x)=f(x+1)$ and
$f(x)+f(x^2+a)=1$ for all $x\in F_q$. An obvious choice is $f=\tr$. In fact,
this is the only natural choice which does not depend on $a$ itself; for we may
assume that $f(0)=0$, and then we must have $f(a)=1$ for all $a$ to which the
discussion applies, namely those $a$ for which $\tr(a)=1$ and $\alpha$ has
order $q+1$. This is close to requiring $f(a)=1$ whenever $\tr(a)=1$, which in
turn implies $f=\tr$ because $f$ must be zero on exactly half the elements
of~$F_q$.

So by the process described we arrive at the definition of the graph $G_k(a)$.

\subsection{Other possibilities.}

Aiming for a circulant is not {\em a priori} the right thing to do; the
Paley graphs are circulants if $q$ is prime but not in general. However, in
order that the group $PSL(2,q)$ have a nice action on our graph we should
choose its edge set to be a union of orbits of elements of $PSL(2,q)$. We
also want the graph to be self-complementary and to have some M\"obius map
interchanging the graph and its complement. As we shall see below, a great
number of M\"obius maps have order $q+1$, and so any graph of this more
general kind is likely to be circulant.

\section{Field Work}\label{field}

Here we make further, standard and elementary, calculations over finite
fields to justify some earlier remarks. A full treatment of these matters
can be found in Lidl and Niederreiter~\cite{LN}.

\subsection{Trace comments.}\label{tracing}
The trace map is defined by $\tr(a)=a+a^2+a^4+\ldots+a^{q/2}$. Thus
$\tr(a)^2=\tr(a)$ so $\tr(a)\in F_2$. Moreover trace is a linear map.
There is a distinction between even $k$ and odd $k$, because
$$ \tr(1)=1+1^2+1^4+\ldots+1^{k-1}=\left\{ \begin{array}{ll}
						0 &\mbox{if $k$ is even} \\
						1 &\mbox{if $k$ is odd.} \\
					    \end{array} \right. $$
Since trace is a linear map,
$$\tr\rat{x}{y}{a}+\tr\rat{y}{x}{a}=\tr(1).$$ 
It follows that the definition in~\S\ref{defn} determines a graph if $k$ is
even and a tournament if $k$ is odd, as claimed.

The map $\tr:F_q\to F_2$ is surjective, since trace, being a polynomial of
degree lower than $q$, cannot annihilate $F_q$. Let $T_i=\tr^{-1}(i)$,
$i=0,1$. Then $T_0$ is the kernel of trace; since the map is surjective, we
have $\dim T_0=k-1$ and so $|T_0|=|T_1|=2^{k-1}=q/2$.

Now $\tr(a^2)=\tr(a)$, or $\tr(a^2+a)=0$.  The map $x\mapsto x^2+x$ is also a
linear map $F_q\to F_q$. Its kernel is $F_2$ so its image has dimension
$k-1$. But its image contains $T_0$. Therefore its image is $T_0$; in
particular, for every element $c$ with $\tr(c)=0$ there exists an element $b\in
F_q$ with $b^2+b=c$. There are two solutions to this quadratic equation, the
other being $b+1$. Thus if $k$ is even and $\tr(1)=0$ the two solutions have
the same trace, whereas if $k$ is odd the solutions have different traces.

\subsection{M\"obius comments.}\label{conjugal}

Our aim here is to identify a suitable element $a\in F_q$ with which to
carry out the above construction. Note that, for any $a$ with $\tr(a)=1$,
then the equation $z^2+z+a=0$ has no solution in $F_q$, because
$\tr(z^2+z+a)=\tr(a)=1$. Therefore the equation has a root $\lambda$ in
$F_{q^2}$. It follows that $\overline\lambda=\lambda^q$ is the other root,
because $\overline\lambda^2 + \overline\lambda + a = (\lambda^2+\lambda +
a)^q$.

Let $k$ be the order of the element $\overline\lambda/\lambda$ in
$F_{q^2}$. Then $1=(\overline\lambda/\lambda)^k=\lambda^{k(q-1)}$, but also
$\lambda^{q^2-1}=1$, so $k\mid (q+1)$ (in particular, if $q+1$ is a Fermat prime
then $\overline\lambda/\lambda$ has order $q+1$).

Now let $k>2$ be any factor of $q+1$. Let $g$ be a primitive root for
$F_{q^2}$. Then the cyclic group $\langle g^{q-1}\rangle$ of order $q+1$ has
exactly $\phi(k)$ elements of order $k$, where $\phi$ is Euler's
function. Let $\mu=g^{t(q-1)}$ be an element of order $k$ in $\langle
g^{q-1}\rangle$. Then $\mu^{q-1}=g^{t(q-1)^2}=g^{-2t(q-1)}$. Therefore
$\mu\notin F_q$, for otherwise $\mu^{q-1}=1$ which would imply $(q+1)\mid 2t$,
which in turn would imply that $\mu^2=1$, contradicting $k>2$.

Given $\mu=g^{t(q-1)}$ as described, let $\nu=g^t$. Let
$b=\nu+\overline\nu$, where $\overline\nu=\nu^q$. Then $\overline b =
b^q=\overline\nu+\nu=b$, so $b\in F_q$. Put $\lambda=\nu/b$. Thus
$\lambda+\overline\lambda=1$, and the element
$\overline\lambda/\lambda=\overline\nu/\nu=\mu$ has order $k$. Let
$\lambda\overline\lambda=a$; since $a^q=a$ we have $a\in F_q$. Moreover
$\lambda^2+\lambda+a=0$, so $\tr(a)=\lambda+\lambda^q=1$.

We summarize as follows. Every element $a$ of trace~1 in $F_q$ satisfies an
equation $\lambda^2+\lambda+a=0$ where $\lambda\in F_{q^2}$ and the order of
$\overline\lambda/\lambda$ divides $q+1$. Conversely, for every factor $k>2$ of
$q+1$ there exists such an $a$ such that $\overline\lambda/\lambda$ has order
$k$. 

In particular, there exists an $a$ such that $\overline\lambda/\lambda$ has
order $q+1$. For such an $a$, consider the map $\alpha:z\to a/(z+1)$ and its
associated matrix $\mat{0}{a}{1}{1}$. This matrix has eigenvectors
$\col\lambda$ and $\col{\overline\lambda}$ with eigenvalues $\overline\lambda$
and $\lambda$ respectively.
Now $\infty=\colinf=\col\lambda+\col{\overline\lambda}$. Therefore the result
of applying the map $z\mapsto a/(z+1)$ to $\infty$ $k$ times is
$\overline\lambda^k\col\lambda+\lambda^k\col{\overline\lambda}$. This equals
$\colinf$ only if $\overline\lambda^k+\lambda^k=0$, which is to say
$(\overline\lambda/\lambda)^k=1$. Since $\overline\lambda/\lambda$ has order
$q+1$, the vertex $\infty$ is in an $\alpha$-orbit of size~$q+1$. Thus there do
exist elements $a$ for which the graph $G_k(a)$ is a self-complementary
circulant graph, as described in~\S\ref{hist}.

\section{Elementary properties}\label{elem}

Some of the more accessible properties of $G_k$ can now be described.

\subsection{Isomorphisms.}\label{isoms}

Let $b\in F_q$. The map $x\mapsto x+b$ is a permutation of $F_q$. Moreover
$$\tr\rat{(x+b)}{(y+b)}{a} = \tr\rat{x}{y}{b^2+b+a}+\tr(b).\eqno{(\dag)}$$ 

Suppose that $k$ is odd, that is, $\tr(1)=1$. Let $a$ and $a^\prime$ be two
elements of $T_1$. Then $\tr(a+a')=0$, and by the remarks
in~\S\ref{tracing}, there exists an element $b$ with $b^2+b+a=a^\prime$ and
$\tr(b)=0$. Therefore $(\dag)$ shows that the map $x\mapsto x+b$ is an
isomorphism $G_k(a^\prime)\to G_k(a)$. Moreover, since $1^2+1+a=a$,
by~$(\dag)$ the map $x\mapsto x+1$ is an orientation-reversing bijection of
the vertex set, because $\tr(1)=1$. It follows that the tournaments defined
in~\S\ref{defn} are isomorphic to each other and are self-complementary.

Now let $k$ be even. We showed in~\S\ref{field} that there is some element
$a$ for which the map $\alpha:z\to a/(z+1)$ has order $q+1$ and
$\tr(a)=1$. Let $a^\prime$ be any other element of $T_1$. Let
$c=a^\prime-a$. Again, by the remarks in~\S\ref{tracing}, there exists an
element $b$ with $b^2+b+a=a'$. Now either $\tr(b)=0$, in which case
$(\dag)$ shows that the map $x\mapsto x+b$ is an isomorphism
$G_k(a^\prime)\to G_k(a)$, or $\tr(b)=1$, in which case the map $x\mapsto
x+b$ is an isomorphism between $G_k(a^\prime)$ and the complement of
$G_k(a)$. But $G_k(a)$ is vertex-transitive and self-complementary, as
shown in~\S\ref{hist}. Therefore the graphs defined in~\S\ref{defn} are
isomorphic to each other, being both vertex-transitive and
self-complementary.

\subsection{Automorphisms.}\label{automs}

Let $a\in F_q$ have trace~one. The M\"obius map $z\mapsto a/(z+1)$ is a
permutation of $V$. It is also an automorphism of the graph $G_k(a)$, because
$$
\frac{\frac{a}{x+1}\cdot\frac{a}{y+1}+\frac{a}{x+1}+a}
{\frac{a}{x+1}+\frac{a}{y+1}} =\frac{xy+x+a}{x+y}.
$$
This, of course, is just the automorphism $\alpha$ that was built
into the definition of $G_k(a)$.

In the graph case, the map $z\mapsto z+1$ is also an automorphism, being
the map $v_i\mapsto v_{-i}$.

\subsection{Co-degrees.}\label{codeg}

The co-degree of a pair $x,y$ of vertices is the number of their common
neighbours. As mentioned earlier, a $q/2$-regular graph of order $q+1$ is
a conference graph if every pair $x,y$ has codegree $q/4-\epsilon$, where
$\epsilon=0$ or~$1$ according as $x$ and $y$ are not adjacent or are
adjacent.

The present graphs do not quite satisfy this condition but come close. Let
us compute the co-degree of $x,y$ in $G_a(q)$. By the rotational symmetry
we may assume that $y=\infty$. A vertex $w\notin\{\infty,x\}$ is joined to
$\infty$ if $\tr(w)=0$ and to $x$ if $\tr((xw+x+a)/(x+w))=0$. Let
$\psi:F_q\to\{-1,1\}$ be the additive character $\psi(z)=(-1)^{\tr(z)}$. If
$\ell$ is the co-degree of $x,y$, then there are $q/2-\epsilon -\ell$
vertices joined to $x$ but not to~$y$, with the same number joined to $y$
but not~$x$. So we have
$$
\sum_{w\in F_q, w\ne x} \psi(w)\psi\rat{x}{w}{a} = q-1
-4(q/2-\epsilon-\ell)\,.
$$
Thus, writing $K$ for the sum on the left, we have
$\ell=q/4-\epsilon+(K+1)/4$.

Using the substitutions $w=z+x$ and $b=x^2+x+a$ we have
$$
K = \sum_{z\in F_q-\{0\}} \psi(z+x)\psi\left(x+\frac{b}{z}\right)
= \sum_{z\in F_q-\{0\}} \psi\left(z+\frac{b}{z}\right)\,.
$$
Therefore $K$ is a {\em Kloosterman sum}; see Lidl and
Niederreiter~\cite[Section~5.5]{LN} for a discussion. In particular,
$|K|\le 2\sqrt q$ (\cite[Theorem~5.45]{LN}). This was proved by Carlitz and
Uchiyama~\cite{CU}, extending the proof by Weil~\cite{W} to even~$q$. A
self-contained proof, based on Stepanov~\cite{St}, appears in
Schmidt~\cite[Chapter~2]{Sc}.

In the case that $G_k$ is a graph we have $q=2^k$ where $k$ is even, and
so $\sqrt q$ is an even integer. Therefore every co-degree is at most
$q/4+\sqrt q / 2$.

\subsection{Pseudo-randomness}

For our purposes, the import of the preceeding estimate of co-degrees is
that the graph $G_k$ is pseudo-random. Specifically, it is $(1/2,
q^{3/4})$-{\em jumbled}, meaning that, for every induced subgraph $H\subset
G_k$, $|e(H)-\frac{1}{2}\binom{|H|}{2}|\le q^{3/4}|H|$ holds. This follows
comfortably from~\cite[Theorem~1.1]{T} using the bound $q/4+\sqrt q/2$ for
co-degrees.

From this it follows that $G_k$ enjoys all the usual consequences of
pseudo-randomness, such as expansion, having about the expected number of
induced subgraphs of any given kind, and so on.

Another approach to showing that $G_k$ is pseudo-random would be to
estimate the eigenvalues, which are of course available in a reasonably
explicit form given that $G_k$ is a circulant. However the present approach
via co-degrees is quick and effective.

\subsection{Hamiltonian decompositions}
As mentioned above, the Paley graph of order $q$ has a Hamiltonian
decomposition when $q$ is prime because it is a circulant of prime order,
and likewise so is $G_k$ if $q+1$ is a Fermat prime, though there seems to
be a limited supply of these.

What about non-prime orders? For the Paley graph, there is always a
Hamiltonian decomposition, as shown by Alspach, Bryant and
Dyer~\cite{ABD}. The graphs $G_k$ too have a Hamiltonian decomposition, at
least if $k$ is large. This follows from the deep work of K\"uhn and
Osthus~\cite{KO}. Theorem~1.2 of~\cite{KO} states that there is some number
$\tau>0$ such that, provided $G_k$ is a {\em robust
  $(\tau/3,\tau)$-expander}, then $G_k$ has a Hamiltonian decomposition for
large~$k$. This condition requires that, for every subset $S$ of the
vertices of~$G_k$ with $\tau q \le |S|\le (1-\tau)q$, there are at least
$|S|+\tau q/3$ vertices of $G_k$ having at least $\tau q/3$ neighbours
in~$S$. The condition is comfortably satisfied by $G_k$ because it is
$(1/2,q^{3/4})$-jumbled (using simple standard properties of such
graphs~\cite{T}). The decomposition is, of course, not explicit but there
is a polynomial time algorithm for finding it.

\section{Acknowledgements}

Thanks are due to Robin Chapman for his comments. In particular he suggests
another description of the graph~$G_k$, from a field theoretic, rather than a
geometric, viewpoint. The line $PG(1,q)$ can be identified in a natural way
with the quotient group $F_{q^2}^*/F_q^*$, so we can consider graphs with
this as its vertex set. Let $\lambda\in F_{q^2}^*$. Given $u,v\in
F_{q^2}^*$ then the equivalence classes $[u], [v]$ are vertices of a graph
$H_\lambda$, and we join $[u]$ to $[v]$ if $\tr(T(\lambda
u^qv)/T(\lambda)T(u^qv))=0$, where $T(x)=x+x^q$ is the trace map
$F_{q^2}^*\to F_q^*$. Chapman~\cite{C} shows that $H_\lambda$ is
well defined and isomorphic to $G_k$.

\end{document}